\documentclass[]{article}

\usepackage{amssymb,latexsym,amsmath,graphicx,subfigure,color,amsthm,IEEEtrantools}     
\usepackage{stmaryrd}
\usepackage{multirow}

\DeclareMathOperator*{\argmin}{argmin}

\newtheorem{theorem}{Theorem}
\newtheorem{lemma}[theorem]{Lemma}

\begin{document}

\title{Constraint Energy Minimizing Generalized Multiscale Discontinuous Galerkin Method}
\author{
Siu Wun Cheung\thanks{Department of Mathematics, Texas A\&M University, College Station, TX 77843, USA (\texttt{tonycsw2905@math.tamu.edu})}
\and
Eric T. Chung\thanks{Department of Mathematics, The Chinese University of Hong Kong, Shatin, New Territories, Hong Kong SAR, China (\texttt{tschung@math.cuhk.edu.hk})}
\and
Wing Tat Leung\thanks{Institute for Computational Engineering and Sciences, 
The University of Texas at Austin, Austin, Texas, USA (\texttt{wleungo@ices.utexas.edu})}
}
\maketitle

\begin{abstract}
Numerical simulation of flow problems and wave propagation in heterogeneous media 
has important applications in many engineering areas. 
However, numerical solutions on the fine grid are often prohibitively expensive, 
and multiscale model reduction techniques are introduced to efficiently solve for 
an accurate approximation on the coarse grid.
In this paper, we propose an energy minimization based multiscale model reduction approach 
in the discontinuous Galerkin discretization setting.  
The main idea of the method is to extract the non-decaying component in the high conductivity regions by 
identifying dominant modes with small eigenvalues of local spectral problems, 
and define multiscale basis functions in coarse oversampled regions by constraint energy minimization problems. 
The multiscale basis functions are in general discontinuous on the coarse grid 
and coupled by interior penalty discontinuous Galerkin formulation. 
The minimal degree of freedom in representing high-contrast features 
is achieved through the design of local spectral problems, 
which provides the most compressed local multiscale space. 
We analyze the method for solving Darcy flow problem 
and show that the convergence is linear in coarse mesh size 
and independent of the contrast, provided that the oversampling size is appropriately chosen. 
Numerical results are presented to show the performance of the method for 
simulation on flow problem and wave propagation in high-contrast heterogeneous media.
\end{abstract}

\section{Introduction}
\label{sec:intro}

Many engineering applications require numerical simulation in heterogeneous media with multiple scales and high contrast.  
For example, Darcy flow equation in heterogeneous media is used to describe 
fluid flow in porous medium in reservoir simulation, 
and wave equation in heterogeneous media has been widely used for subsurface modeling. 
Numerical solutions on the fine grid are often prohibitively expensive in these complex multiscale problems. 

To this end, extensive research effort had been devoted to developing efficient methods 
for solving multiscale problems at reduced expense, 
for example, numerical homogenization approaches \cite{papanicolau1978asymptotic,weh02}
and multiscale methods, including
Multiscale Finite Element Methods (MsFEM) \cite{hw97,ehw99,ch03,eh09,cgh09}, 
Variational Multiscale Methods (VMS) \cite{hfmq98,hughes2007variational,Iliev_MMS_11,calo2011note}, 
Heterogeneous Multiscale Methods (HMM) \cite{ee03,abdulle05,emz05} and 
and Generalized Multsicale Finite Element Methods (GMsFEM) 
\cite{egh12,chung2014adaptive,chung2015generalizedwave,chung2016adaptive}. 
The common goal of these methods is to construct numerical solvers 
on the coarse grid, which is typically much coarser than the fine grid
which captures all the heterogeneities in the medium properties. 
In numerical homogenization approaches, effective properties are computed and 
the global problem is formulated and solved on the coarse grid. 
However, these approaches are limited to the cases when the medium properties possess scale separation.
On the other hand, multiscale methods construct of multiscale basis functions which are 
responsible for capturing the local oscillatory effects of the solution. 
Once the multiscale basis functions, coarse-scale equations are formulated.
Moreover, fine-scale information can be recovered by the coarse-scale coefficients and mutliscale basis functions. 

Many existing mutliscale methods, such as MsFEM, VMS and HMM, 
construct one basis function per local coarse region to handle the effects of local heterogeneities. 
However, for more complex multiscale problems, each local coarse region contains several high-conductivity regions and 
multiple multiscale basis functions are required to represent the local solution space. 
GMsFEM is developed to allow systematic enrichment of the coarse-scale space with fine-scale information 
and identify the underlying low-dimensional local structures for solution representation. 
The main idea of GMsFEM is to extract local dominant modes by carefully designed local spectral problems in coarse regions, 
and the convergence of the GMsFEM is related to eigenvalue decay of local spectral problems. 
For a more detailed discussion on GMsFEM, we refer the readers to 
\cite{egw10,egh12,eglp13oversampling,chung2014adaptive,chung2015residual,
chung2016adaptive,efendiev2017bayes,cheung2019bayes,park2019mc,wang2020vug} 
and the references therein. 
Our method developed in this work is motivated by GMsFEM and achieves spectral convergence. 
Through the design of local spectral problems, our method results in the minimal degree of freedom 
in representing high-contrast features. 

For typical mesh-based numerical discretization, such as the finite element method and finite difference method,
the most important issue of the solution accuracy is mesh convergence. 
However, for multiscale problems, it is difficult to adjust coarse-grid size based on scales and contrast, 
and it is desirable to have convergence independent of these physical parameters. 
it becomes non-trivial to derive multiscale methods with convergence on coarse mesh size independent of scales and contrast.  
Recently, several multiscale methods with mesh convergence are developed. 
\cite{owhadi2014polyharmonic, maalqvist2014localization, owhadi2017multigrid}. 
This idea can be combined with the use of local spectral problems for 
achieving both spectral convergence and mesh convergence 
\cite{hou2017sparse,chung2018constraint,chung2018mixed,cheung2018mc}. 
On the other hand, to overcome the difficulty of stability and conservation for 
convection-dominated problems and wave propagations in heterogeneous media, 
multiscale methods in the discontinuous Galerkin (DG) framework have been investigated 
\cite{ehw99,buffa2006analysis,riviere2008discontinuous,eglmsMSDG,elfverson2013dg,efendiev2015spectral,chung2017dg,chung2018dg}. 
In these approaches, unlike conforming finite element formulations, 
multiscale basis functions are in general discontinuous on the coarse grid, 
and stabilization or penalty terms are added to ensure well-posedness of the global problem. 
Our goal in this work is to develop a robust multiscale method in interior penalty discontinuous Galerkin formulation, 
which exhibits both spectral convergence and mesh convergence. 

In this paper, we present Constraint Energy Minimizing Generalized Multiscale
Discontinuous Galerkin Method (CEM-GMsDGM). 
There are two key ingredients of the presented approach.  
The first main ingredient is the local spectral problems in each coarse block for identification of auxiliary basis functions. 
The low-energy dominant modes, which are eigenvectors corresponding to small eigenvalues of local spectral problems, 
are used as auxiliary basis functions for further construction. 
The auxiliary basis functions possess the information related to high conductivity channels and 
it suffices to use the same number of auxiliary basis functions as the number of channels in a coarse block.
The second ingredient is the constraint energy minimization problems for definition of multiscale basis functions. 
Each of the auxiliary basis functions sets up an independent constraint and 
uniquely defines a corresponding multiscale basis function. 
The multiscale basis functions will then be used to span the multiscale space and 
used to solve the coarse-scale global problem in IPDG formulation.
We remark that the local spectral problems and the constraint energy minimization problems 
are carefully designed and supported by our analysis. 
Thanks to the design of local spectral problems, 
the auxiliary space is of minimal dimension for representing high-contrast features 
and obtaining a contrast-independent convergence. 
Due to the fact that the dimensions of the auxiliary space and the multiscale space are identical, 
the multiscale space is of minimum dimension as well. 
In the construction of multiscale basis functions, the constraints are responsible for 
handling non-decaying components represented by the auxiliary basis functions in the high conductivity regions 
and achieving linear convergence in coarse mesh size. 
On the other hand, the multiscale basis functions are supported in oversampled coarse regions 
and allowed to have discontinuity on the coarse grid. 
Therefore, the IPDG bilinear form is also used to define the energy term in the constraint energy minimization problems. 
The advantages of the method is verified both theoretically and numerically. 
We analyze the method for solving Darcy flow problem 
and establish a criterion for the oversampling size 
which is sufficient for linear coarse-mesh convergence independent of the contrast.
Numerical results are presented to show the performance of the method for 
simulation on flow problem and wave propagation in high-contrast heterogeneous media. 

The paper is organized as follows. In Section~\ref{sec:prelim}, we will introduce the notions of grids, 
and essential discretization details such as DG finite element spaces and IPDG formulation on the coarse grid. 
The details of the proposed method will be presented in Section~\ref{sec:method}. 
The method will be analyzed in Section~\ref{sec:analysis}.
Numerical results will be provided in Section~\ref{sec:numerical}.
Finally, a conclusion will be given in Section~\ref{sec:conclusion}.

\section{Preliminaries}
\label{sec:prelim}
We consider the following high-contrast flow problem
\begin{equation}
- \text{div} \left(\kappa \nabla u\right) = f \text{ in } \Omega,
\label{eq:elliptic}
\end{equation}
subject to the homogeneous Dirichlet boundary condition $u = 0$ on $\partial \Omega$, 
where $\Omega \subset \mathbb{R}^2$ is the computational domain and $f$ is a given source term. 
We assume that the permeability field $\kappa$ is highly heterogeneous 
with very high contrast $\kappa_0 \leq \kappa \leq \kappa_1$. 

\begin{figure}[ht!]
\centering
\includegraphics[width=0.5\linewidth]{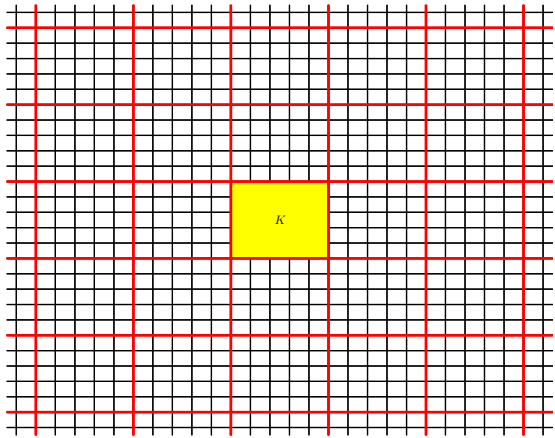}
\caption{An illustration of the fine grid and the coarse grid and a coarse element.}
\label{fig:mesh}
\end{figure}

Next, we introduce the notions of coarse and fine meshes. 
We start with a usual partition $\mathcal{T}^H$ of $\Omega$ into finite elements, 
which does not necessarily resolve any multiscale features. 
The partition $\mathcal{T}^H$ is called a coarse grid and 
a generic element $K$ in the partition $\mathcal{T}^H$ is called a coarse element. 
Moreover, $H > 0$ is called the coarse mesh size.
We let $N_c$ be the number of coarse grid nodes and 
$N$ be the number of coarse elements. 
We also denote the collection of all coarse grid edges by $\mathcal{E}^H$.
We perform a refinement of $\mathcal{T}^H$ to obtain a fine grid $\mathcal{T}^h$, 
where $h > 0$ is called the fine mesh size. 
It is assumed that the fine grid is sufficiently fine to resolve the solution. 
An illustration of the fine grid and the coarse grid and a coarse element are shown in Figure~\ref{fig:mesh}. 

We are now going to discuss the discontinuous Galerkin (DG) discretization 
and the interior penalty discontinuous Galerkin (IPDG) global formulation. 
For the $i$-th coarse block $K_i$, we denote the restriction of the Sobolev space $H_0^1(\Omega)$
on $K_i$ by $V(K_i)$. We let $V_h(K_i)$ be 
the conforming bilinear elements defined on the fine grid $\mathcal{T}^h$ in $K_i$, i.e.
\begin{equation}
V_h(K_i) = \left\{ v \in V(K_i): v \vert_\tau \in \mathbb{Q}^1(\tau) \text{ for all } 
\tau \in \mathcal{T}^h \text{ and } \tau \subset K_i\right\}, 
\end{equation}
where $\mathbb{Q}^1(\tau)$ stands for the bilinear element on the fine grid block $\tau$. 
The DG approximation space is then given by the space of 
coarse-scale locally conforming piecewise bilinear 
fine-grid basis functions, namely 
\begin{equation}
V_h = \oplus_{i=1}^N V_h(K_i). 
\end{equation}
We remark that functions in $V_h$ are continuous within coarse blocks, 
but discontinuous across the coarse grid edges in general. 
The global formulation of IPDG method then reads:
find $u_h \in V_h$ such that
\begin{equation}
a_{{DG}}\left(u_h ,w\right) = \int_{\Omega} fw \, dx \text{ for all } w \in V_h, 
\label{eq:sol_dg}
\end{equation}
where the bilinear form $a_{{DG}}$ is defined by:
\begin{equation}
\begin{split}
a_{{DG}}\left(v,w\right) 
& = \sum_{K \in \mathcal{T}^H} \int_K \kappa \nabla v \cdot \nabla w \, dx 
- \sum_{E \in \mathcal{E}^H} \int_E \{ \kappa \nabla v \cdot n_E \} \llbracket w \rrbracket \, d\sigma \\
& \quad - \sum_{E \in \mathcal{E}^H} \int_E \{ \kappa \nabla w \cdot n_E \} \llbracket v \rrbracket \, d\sigma 
+ \dfrac{\gamma}{h} \sum_{E \in \mathcal{E}^H} \int_E \overline{\kappa} \llbracket v \rrbracket \llbracket w \rrbracket \, d\sigma,
\end{split}
\label{eq:dg_bilinear}
\end{equation}
where $\gamma > 0$ is a penalty parameter and 
$n_E$ is a fixed unit normal vector defined on the coarse edge $E \in \mathcal{E}^H$. 
Note that, in \eqref{eq:dg_bilinear}, the average and the jump operators 
are defined in the classical way. 
Specifically, consider an interior coarse edge $E \in \mathcal{E}^H$ and 
let $K^+$ and $K^-$ be the two coarse grid blocks sharing the edge $E$, 
where the unit normal vector $n_E$ is pointing from $K^+$ to $K^-$. 
For a piecewise smooth function $G$ with respect to the coarse grid $\mathcal{T}^H$, we define
\begin{equation}
\begin{split}
\{ G \} & = \dfrac{1}{2}\left(G^+ + G^-\right), \\
\llbracket G \rrbracket & = G^+ - G^-,
\end{split}
\end{equation}
where $G^+ = G\vert_{K^+}$ and $G^- = G\vert_{K^-}$. 
Moreover, on the edge $E$, we define $\overline{\kappa} = \left(\kappa_{K^+} + \kappa_{K^-} \right)/2$, 
where $\kappa_{K^\pm}$ is the maximum value of $\kappa$ over $K^\pm$. 
For a coarse edge $E$ lying on the boundary $\partial \Omega$, we define
$\{G\}=\llbracket G \rrbracket =G$, and $\kappa=\kappa_K$ on $E$, 
where we always assume that $n_E$ is pointing outside of $\Omega$. 

First, we define the energy norm on the space $V$ of coarse-grid piecewise smooth functions by
\begin{equation}
\| w \|_a^2 = a_{DG}(w,w) \text{ for all } w \in V. 
\end{equation}
We also define the DG-norm on $V$ by 
\begin{equation}
\| w \|_{DG}^2 = \sum_{K \in \mathcal{T}^H} \int_K \kappa \vert \nabla w \vert^2 \, dx 
+ \dfrac{\gamma}{h} \sum_{E \in \mathcal{E}^H} \int_E \overline{\kappa} \llbracket w \rrbracket^2 \, d\sigma 
\text{ for all } w \in V.
\end{equation}
The two norms are equivalent on the subspace of piecewise bi-cubic polynomials in $V$: 
there exists $C_0 \geq 1$ such that 
\begin{equation}
C_0^{-1} \| w \|_{a} \leq \| w \|_{DG} \leq C_0 \| w \|_a.
\end{equation}
The continuity and coercivity results of the bilinear form $a_{DG}$ 
with respect to the DG-norm is ensured by a sufficiently large penalty parameter $\gamma$.
While the method works well for general highly heterogeneous field $\kappa$, 
we assume $\kappa$ is piecewise constant on the fine grid $\mathcal{T}^h$ 
for the sake of simplicity in our analysis presented in Section~\ref{sec:analysis}. 

\section{Method description}
\label{sec:method}
\begin{figure}[ht!]
\centering
\includegraphics[width=0.5\linewidth]{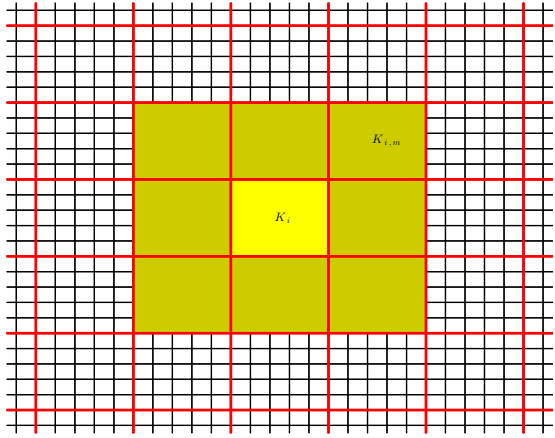}
\caption{An illustration of an oversampled domain formed by enlarging $K_i$ with $1$ coarse grid layer.}
\label{fig:oversample}
\end{figure}

In this section, we will present the construction of the multiscale basis functions. 
First, we will use the concept of GMsFEM to construct our auxiliary multiscale basis functions 
on a generic coarse block $K$ in the coarse grid. 
We consider $V_h(K_i)$ as the snapshot space in $K_i$ and 
perform a dimension reduction through a spectral problem, 
which is to find a real number $\lambda_j^{\left(i\right)}$ 
and a function $\phi_j^{\left(i\right)} \in V_h(K_i)$ such that
\begin{equation}
a_i\left(\phi_j^{\left(i\right)}, w\right) = \lambda_j^{\left(i\right)} s_i\left(\phi_j^{\left(i\right)}, w\right) \text{ for all } w \in V_h(K_i),
\label{eq:spectral_prob}
\end{equation}
where $a_i$ is a symmetric non-negative definite bilinear operator 
and $s_i$ is a symmetric positive definite bilinear operators defined on $V_h(K_i) \times V_h(K_i)$. 
We remark that the above problem is solved on the fine mesh in the actual computations. 
Based on our analysis, we can choose
\begin{equation}
\begin{split}
a_i\left(v,w\right) & = \int_{K_i} \kappa \nabla v \cdot \nabla w \, dx, \\
s_i\left(v,w\right) & = \int_{K_i} \widetilde{\kappa} v w \, dx,
\end{split}
\label{eq:spectral_bilinear_form}
\end{equation}
where $\tilde{\kappa} = \sum_{j=1}^{N_c} \kappa \vert \nabla \chi_j^{ms} \vert^2$ and 
$\{\chi_j^{ms}\}_{j=1}^{N_c}$ are the standard multiscale finite element (MsFEM) basis functions.
We let $\lambda_j^{\left(i\right)}$ be the eigenvalues of \eqref{eq:spectral_prob} 
arranged in ascending order in $j$, and use the first $L_i$ eigenfunctions 
to construct our local auxiliary multiscale space
\begin{equation}
V_{aux}^{\left(i\right)} = \text{span} \{ \phi_j^{\left(i\right)}: 1 \leq j \leq L_i\}.
\end{equation}
The global auxiliary multiscale space $V^h_{aux}$ is then defined as 
the sum of these local auxiliary multiscale spaces
\begin{equation}
V_{aux} = \oplus_{i=1}^N V_{aux}^{\left(i\right)}.
\end{equation}
For the local auxiliary multiscale space $V^{\left(i\right)}_{aux}$, the bilinear form $s_i$ 
in \eqref{eq:spectral_bilinear_form} defines an inner product 
with norm $\|v\|_{s(K_i)} = s\left(v,v\right)^\frac{1}{2}$. 
These local inner products and norms provide a natural definitions of inner
product and norm for the global auxiliary multiscale space $V_{aux}$, which are defined by
\begin{equation}
\begin{split}
s\left(v,w\right) & = \sum_{i=1}^N s_i\left(v,w\right) \text{ for all } v,w \in V_{aux}, \\
\|v\|_s & = s\left(v,v\right)^\frac{1}{2} \text{ for all } v \in V_{aux}.
\end{split}
\end{equation}
We note that $s\left(v, w\right)$ and $\|v\|_s$ are also an inner product and norm for the space $V_h$.
Before we move on to discuss the construction of multiscale basis functions, 
we introduce some tools which will be used to describe our method and analyze the convergence.
We first introduce the concept of $\phi$-orthogonality. 
For $1 \leq i \leq N$ and $1 \leq j \leq L_{i}$, 
in coarse block $K_i$, given auxiliary basis function $\phi_j^{\left(i\right)} \in V_{aux}$,
we say that $\psi \in V_h$ is $\phi_j^{\left(i\right)}$-orthogonal if
\begin{equation}
s \left(\psi, \phi_{j'}^{\left(i'\right)}\right) = \delta_{i,i'} \delta_{j,j'} \text{ for all } 1 \leq j' \leq L_{i'} \text{ and } 1 \leq i' \leq N.
\end{equation}
We also introduce a projection operator $\pi: V_h \to V_{aux}$ by 
$\pi = \sum_{i=1}^N \pi_i$, where
\begin{equation}
\pi_i(v) = \sum_{j=1}^{L_i} \dfrac{s_i\left(v,\phi_j^{\left(i\right)}\right)}{s_i\left(\phi_j^{\left(i\right)},\phi_j^{\left(i\right)}\right)} \phi_j^{\left(i\right)} \text{ for all } v \in V_h, \text{ for all } i = 1,2,\ldots,N.
\end{equation}
Next, we construct our global multiscale basis functions in $V_h$.
The global multiscale basis function $\psi_{j}^{\left(i\right)} \in V_h$ is defined as the solution of 
the following constrained energy minimization problem
\begin{equation}
\psi_{j}^{\left(i\right)} = \argmin \left\{ a_{{DG}}\left(\psi, \psi\right) : 
\psi \in V_h \text{ is } \phi_j^{\left(i\right)} \text{-orthogonal}\right\}.
\label{eq:min1_glo}
\end{equation}
By introducing a Lagrange multiplier, 
the minimization problem \eqref{eq:min1_glo} is equivalent to the following variational problem: 
find $\psi_{j}^{\left(i\right)} \in V_h$ and $\mu_{j}^{\left(i\right)} \in V_{aux}^{\left(i\right)}$ such that
\begin{equation}
\begin{split}
a_{{DG}}\left(\psi_{j}^{\left(i\right)}, \psi\right) + s_i\left(\psi, \mu_{j}^{\left(i\right)}\right) & = 0 \text{ for all } \psi \in V_h, \\
s_i\left(\psi_{j}^{\left(i\right)} - \phi_j^{\left(i\right)}, \mu\right) & = 0 \text{ for all } \mu \in V_{aux}^{\left(i\right)}.
\end{split}
\label{eq:var1_glo}
\end{equation}
Now we discuss the construction our localized multiscale basis functions. 
We first denote by $K_{i,m}$ an oversampled domain formed by 
enlarging the coarse grid block $K_i$ by $m$ coarse grid layers. 
An illustration of an oversampled domain is shown in Figure~\ref{fig:oversample}.
We introduce the subspace $V_h\left(K_{i,m}\right)$, which contains restriction of 
fine-scale basis functions in $V_h$ on the oversampled domain $K_{i,m}$. 
We also define $V_{h,0}\left(K_{i,m}\right) = V_h(K_{i,m}) \cap H_0^1(K_{i,m})$ by 
the subspace of functions in $V_h\left(K_{i,m}\right)$ vanishing on the 
boundary of the oversampled domain $K_{i,m}$. 
Motivated by the construction of our global multiscale basis functions,
the method for construction of the localized multiscale basis functions are as follows:
The localized multiscale basis function $\psi_{j,{ms}}^{\left(i\right)} \in V_h\left(K_{i,m}\right)$ is defined as the solution of 
the following constrained energy minimization problem
\begin{equation}
\psi_{j,{ms}}^{\left(i\right)} = \argmin \left\{ a_{{DG}}\left(\psi, \psi\right) : 
\psi \in V_h\left(K_{i,m}\right) \text{ is } \phi_j^{\left(i\right)} \text{-orthogonal}\right\}.
\label{eq:min1}
\end{equation}
Using the method of Lagrange multiplier, 
the minimization problem \eqref{eq:min1} is equivalent to the following variational problem: 
find $\psi_{j,{ms}}^{\left(i\right)} \in V_h\left(K_{i,m}\right)$ and $\mu_{j}^{\left(i\right)} \in V_{aux}^{\left(i\right)}$ such that
\begin{equation}
\begin{split}
a_{{DG}}\left(\psi_{j,{ms}}^{\left(i\right)}, \psi\right) + s_i\left(\psi, \mu_{j}^{\left(i\right)}\right) & = 0 \text{ for all } \psi \in V_h\left(K_{i,m}\right), \\
s_i\left(\psi_{j,{ms}}^{\left(i\right)} - \phi_j^{\left(i\right)}, \mu\right) & = 0 \text{ for all } \mu \in V_{aux}^{\left(i\right)}.
\end{split}
\label{eq:var1}
\end{equation}
We use the localized multiscale basis functions to construct
the multiscale DG finite element space, which is defined as
\begin{equation}
V_{ms} = \text{span}  \{\psi_{j,{ms}}^{\left(i\right)} : 1 \leq j \leq L_i, 1 \leq i \leq N \}.
\label{eq:msdg_space}
\end{equation}
We remark that the multiscale finite element space $V_{ms}$ is a subspace of $V_h$. 
After the multiscale DG finite element space is constructed, 
the multiscale solution $u_{ms}$ is given by: find $u_{ms} \in V_{ms}$ such that
\begin{equation}
a_{{DG}}\left(u_{ms},w\right) = \int_{\Omega} fw \, dx \text{ for all } w \in V_{ms}. 
\label{eq:sol_ms}
\end{equation}

\section{Analysis}
\label{sec:analysis}
In this section, we will analyze the proposed method.
Besides the energy norm and the DG norm, we also define the $s$-norm on $V$ by
\begin{equation}
\| w \|_s^2 = \sum_{K \in \mathcal{T}^H} \int_K \widetilde{\kappa} \vert w \vert^2 \, dx.
\end{equation}
Given a subdomain $\Omega' \subseteq \Omega$ formed by a union of coarse blocks $K \in \mathcal{T}^H$, 
we also define the local $s$-norm by 
\begin{equation}
\| w \|_{s\left(\Omega'\right)}^2 
= \sum_{K \subseteq \Omega'} \int_K \widetilde{\kappa} \vert w \vert^2 \, dx. 
\end{equation}

The flow of our analysis goes as follows. 
First, we prove the convergence using the global multiscale basis functions.
With the global multiscale basis functions constructed, 
the global multiscale finite element space is defined by
\begin{equation}
V_{glo} = \text{span}  \{\psi_{j}^{\left(i\right)} : 1 \leq j \leq L_i, 1 \leq i \leq N \}, 
\end{equation}
and an approximated solution $u_{glo} \in V_{glo}$ is given by
\begin{equation}
a_{DG}\left(u_{glo},w\right) = \int_{\Omega} fw \, dx \text{ for all } w \in V_{glo}. 
\label{eq:sol_glo}
\end{equation}
We remark that the construction of global multiscale basis functions 
motivates the construction of localized multiscale basis functions.
The approximated solution $u_{glo}$ will also be used in our convergence analysis.
Next, we give an estimate of the difference between 
the global multiscale functions $\psi_j^{\left(i\right)}$ and 
the localized multiscale basis functions $\psi_{j,ms}^{\left(i\right)}$, 
in order to show that using the multiscale solution $u_{ms}$ 
provide similar convergence results as the global solution $u_{glo}$.
For this purpose, we denote the kernel of 
the projection operator $\pi$ by $\widetilde{V}_h$. 
Then, for any $\psi_j^{\left(i\right)} \in V_{glo}$, we have 
\begin{equation}
a_{DG}\left(\psi_j^{\left(i\right)}, w\right) = 0 \text{ for all } w \in \widetilde{V}_h,
\end{equation}
which implies $\widetilde{V}_h \subseteq V_{glo}^\perp$, 
where $V_{glo}^\perp$ is the orthogonal complement of 
$V_{glo}$ with respect to the inner product $a_{DG}\left(\cdot, \cdot\right)$. 
Moreover, since $\text{dim}\left(V_{glo}\right) = \text{dim}\left(V_{aux}\right)$, 
we have $\widetilde{V}_h = V_{glo}^\perp$ and 
$V_h = V_{glo} \oplus \widetilde{V}_h$.

\subsection{Convergence result}

The convergence analysis will start with the following lemma, 
which concerns about the convergence of the approximated solution 
by the global multiscale basis functions.
\begin{lemma}
\label{lemma1}
Let $u_h \in V_h$ be the solution of \eqref{eq:sol_dg} 
and $u_{glo} \in V_{glo}$ be the solution of \eqref{eq:sol_glo}
with the global multiscale basis functions defined by 
the constrained energy minimization problem \eqref{eq:min1_glo}.
Then we have $u_h - u_{glo} \in \widetilde{V}_h$ and
\begin{equation}
\label{eq:lemma1.0}
\| u_h - u_{glo} \|_a \leq \Lambda^{-\frac{1}{2}} \| \widetilde{\kappa}^{-\frac{1}{2}} f \|_{L^2(\Omega)},
\end{equation}
where
\begin{equation}
\Lambda = \min_{1 \leq i \leq N} \lambda_{L_i+1}^{\left(i\right)}.
\end{equation}
Moreover, if we replace the multiscale partition of unity $\{ \chi_j^{ms} \}$ by 
the bilinear partition of unity, we have
\begin{equation}
\| u_h - u_{glo} \|_a \leq C H \Lambda^{-\frac{1}{2}} \| \kappa^{-\frac{1}{2}} f \|_{L^2(\Omega)}.
\end{equation}
\begin{proof}
By the definitions of $u_h$ in \eqref{eq:sol_dg} and $u_{glo}$ in \eqref{eq:sol_glo}, we have
\begin{equation}
\begin{split}
a_{DG}\left(u_h,w\right) & = \int_{\Omega} fw \, dx \text{ for all } w \in V_h, \\
a_{DG}\left(u_{glo},w\right) & = \int_{\Omega} fw \, dx \text{ for all } w \in V_{glo}. 
\end{split}
\label{eq:lemma1.1}
\end{equation}
Since $V_{glo} \subseteq V_h$, this yields the Galerkin orthogonality property.
\begin{equation}
a_{DG}\left(u_h - u_{glo},w\right) = 0 \text{ for all } w \in V_{glo},
\label{eq:lemma1.2}
\end{equation}
which implies $u_h - u_{glo} \in V_{glo}^\perp = \widetilde{V}_h$. 
In particular, if we take $w = u_{glo}$ in \eqref{eq:lemma1.2}, together with \eqref{eq:sol_dg}, we have
\begin{equation}
\begin{split}
\left\| u_h - u_{glo} \right\|_a^2
& = a_{DG}\left(u_h, u_h - u_{glo}\right) \\
& = \left(f, u_h - u_{glo}\right)_{0,\Omega} \\
& \leq \| \widetilde{\kappa}^{-\frac{1}{2}} f \|_{L^2(\Omega)} \| u_h - u_{glo} \|_s.
\end{split}
\label{eq:lemma1.3}
\end{equation}
Since $u_h - u_{glo} \in \widetilde{V}_h$, we have $\pi\left(u_h - u_{glo}\right) = 0$. 
Furthermore, since $K_i$ are disjoint, we have $\pi_i\left(u_h - u_{glo}\right) = 0$ for all $i = 1,2,\ldots,N$.
This implies
\begin{equation}
\begin{split}
\| u_h - u_{glo} \|_s^2
& = \sum_{i=1}^N \| u_h - u_{glo} \|_{s(K_i)}^2 \\
& = \sum_{i=1}^N \| \left(I - \pi_i\right) \left(u_h - u_{glo}\right) \|_{s(K_i)}^2 \\
\end{split}
\label{eq:lemma1.4}
\end{equation}
By the $s_i$-orthogonality of the eigenfunctions $\phi_j^{\left(i\right)}$, we have
\begin{equation}
\begin{split}
\| \left(I - \pi_i\right) \left(u_h - u_{glo}\right) \|_{s(K_i)}^2 
& \leq \left(\lambda^{\left(i\right)}_{L_i+1}\right)^{-1} a_i(u_h - u_{glo},u_h - u_{glo}) \\
& \leq \Lambda^{-1} a_i(u_h - u_{glo},u_h - u_{glo}).
\end{split}
\label{eq:lemma1.5}
\end{equation}
Therefore, we have
\begin{equation}
\begin{split}
\| u_h - u_{glo} \|_s^2
& \leq \Lambda^{-1} \sum_{i=1}^N a_i(u_h - u_{glo},u_h - u_{glo}) \\
& \leq \Lambda^{-1} \| u_h - u_{glo} \|_a^2.
\end{split}
\label{eq:lemma1.6}
\end{equation}
Using \eqref{eq:lemma1.3} and \eqref{eq:lemma1.6}, we obtain our desired result. 
The second part of the result follows from the property 
$\vert \nabla \chi_j \vert = O\left(H^{-1}\right)$ 
of the bilinear partition of unity.
\end{proof}
\end{lemma}

The next step is to prove the global basis functions are indeed localizable. 
This makes use of the following lemma, which states some approximation properties 
of the projection operator $\pi$. 
In the analysis, we will make use of the Lagrange interpolation operator 
and a bubble function in the coarse grid. 
We define the Lagrange interpolation operator
$I_h: C^0(\Omega) \cap H_0^1(\Omega) \to C^0(\Omega) \cap V_h$ by:
for all $u \in C^0(\Omega) \cap H_0^1(\Omega)$, 
the interpolant $I_h u \in C^0(\Omega) \cap V_h$ is defined 
piecewise on each fine block $\tau \in \mathcal{T}^h$ by
\begin{equation}
(I_h u)(x) = u(x) \text{ for all vectices } x \text{ of } \tau,
\end{equation}
which satisfies the standard approximation properties: 
there exists $C_I \geq 1$ such that for every $u \in C^0(\Omega) \cap H_0^1(\Omega)$, 
\begin{equation}
\left\| \widetilde{\kappa}^\frac{1}{2} (u - I_h u) \right\|_{L^2(\tau)} + h \left\| \kappa^\frac{1}{2} \nabla \left(u - I_h u\right) \right\|_{L^2(\tau)} \leq C_{I} h \left\| \kappa^\frac{1}{2} \nabla u \right\|_{L^2(\tau)},
\end{equation}
on each fine block $\tau \in \mathcal{T}^h$. 
For any coarse grid block $K$,  we define a bubble function $B$ on $K$, i.e.
$B(x) = 0$ for all $x \in \partial K$ and $B(x) > 0$ for all $x \in \text{int}\left(K\right)$. 
More precisely, we take $B =  \prod_j \chi_j^{ms}$, where the product is taken 
over all the coarse grid nodes lying on the boundary $\partial K$. 
We can then define the constant 
\begin{equation}
C_\pi = \sup_{K \in \mathcal{T}^H, \mu \in V_{aux}} 
\dfrac{\int_K \widetilde{\kappa} \mu^2}{\int_K \widetilde{\kappa} B \mu^2}.
\end{equation}
In our following analysis, we will assume 
the following smallness criterion on the fine mesh size $h$:
\begin{equation}
C_\pi C_I (C_\mathcal{T}^2 + \lambda_{max}) \| \Theta \|_{L^\infty(\Omega)}^\frac{1}{2} h < 1,
\label{eq:fine-small}
\end{equation}
where $C_\mathcal{T}$ is the maximum number of vertices over all coarse elements $K \in \mathcal{T}^H$ and
\begin{equation}
\begin{split}
\lambda_{max} & = \max_{1 \leq i \leq N} \lambda_{L_i}^{(i)}, \\
\Theta & = \sum_j \vert \nabla \chi_j^{ms} \vert^2. 
\end{split}
\end{equation}

\begin{lemma}
\label{lemma2}
Assume the smallness criterion \eqref{eq:fine-small} on the fine mesh size $h$. 
For any $v_{aux} \in V_{aux}$, there exists a function $v \in C^0(\Omega) \cap V_h$ such that 
\begin{equation}
\pi(v) = v_{aux}, \quad 
\| v \|_a^2 \leq D \| v_{aux} \|_s^2, \quad
\text{supp}(v) \subseteq \text{supp}(v_{aux}),
\end{equation}
where the constant $D$ is defined by 
\begin{equation}
D = \left(\dfrac{2C_\pi(1+C_I^2) \left( C_\mathcal{T}^2 + \lambda_{max} \right)}
{1 - C_\pi C_I \left( C_\mathcal{T}^2+ \lambda_{max} \right) \| \Theta\|_{L^\infty(K_i)}^\frac{1}{2} h}\right)^2.
\end{equation}
\begin{proof}
Let $v_{aux} \in V_{aux}^{(i)}$. 
We consider the following constraint minimization problem on the block $K_i$: 
\begin{equation}
v = \argmin \left\{ a_{DG}(v,v): \, v \in V_{h,0}(K_i), 
\quad s_i(v,\nu) = s_i(v_{aux}, \nu) \text{ for all } \nu \in V_{aux}^{(i)}\right\}.
\label{eq:lemma2.0min}
\end{equation}
The minimization problem \eqref{eq:lemma2.0min} is equivalent to the following variational problem: 
find $(v,\mu) \in V_{h,0}(K_i) \times V_{aux}^{(i)}$ such that 
\begin{equation}
\begin{split}
a_i\left(v, w\right) + s_i\left(w, \mu\right) & = 0 \text{ for all } w \in V_{h,0}(K_i), \\
s_i\left(v - v_{aux}, \nu\right) & = 0 \text{ for all } \nu \in V_{aux}^{\left(i\right)}.
\end{split}
\label{eq:lemma2.0var}
\end{equation}
The existence of solution of \eqref{eq:lemma2.0var} is based on an inf-sup condition: 
\begin{equation}
\inf_{\nu \in V_{aux}^{(i)}} \sup_{w \in V_{h,0}(K_i) \setminus \{0\}} \dfrac{s_i(w,\nu)}{a_i(w,w)} \geq \beta, 
\label{eq:inf-sup}
\end{equation}
where $\beta > 0$ is a constant independent to be determined. 
Pick any $\nu \in V_{aux}^{(i)}$. 
We take $w = I_h(B \nu) \in C^0(\Omega) \cap V_h$. Since $\nu \in V_h(K_i)$ and 
$B(x) = 0$ for all vertices of $K_i$, we have $w \in V_{h,0}(K_i)$. 
First, we see that 
\begin{equation}
\begin{split}
s_i(w,\nu) 
& = \int_{K_i} \widetilde{\kappa} B \nu^2 + \int_{K_i}\widetilde{\kappa}(I_h (B \nu) - B\nu)\nu \\
& \geq C_\pi^{-1} \| \nu \|_{s(K_i)}^2 - \left\| I_h (B \nu) - B\nu \right\|_s \| \nu \|_{s(K_i)} \\
& \geq C_\pi^{-1} \| \nu \|_{s(K_i)}^2 - C_I h \left\| \widetilde{\kappa}^\frac{1}{2} \nabla (B\nu) \right\|_{L^2(K_i)} \| \nu \|_{s(K_i)} \\
& \geq C_\pi^{-1} \| \nu \|_{s(K_i)}^2 - C_I \| \Theta\|_{L^\infty(K_i)}^\frac{1}{2} h
\left\| \kappa^\frac{1}{2} \nabla (B\nu) \right\|_{L^2(K_i)} \| \nu \|_{s(K_i)}.
\end{split}
\label{eq:inf-sup1}
\end{equation}
On the other hand, we observe that 
\begin{equation}
\begin{split}
a_i(w,w) & \leq 2\left(\left\| \kappa^\frac{1}{2} \nabla (B\nu) \right\|_{L^2(K_i)}^2
+ \left\| \kappa^\frac{1}{2} \nabla (B\nu - I_h(B\nu)) \right\|_{L^2(K_i)}^2\right)  \\
& \leq 2(1+C_I^2) \left\| \kappa^\frac{1}{2} \nabla (B\nu) \right\|_{L^2(K_i)}.
\end{split}
\label{eq:inf-sup2}
\end{equation}
It remains to estimate the term $\left\| \kappa^\frac{1}{2} \nabla (B\nu) \right\|_{L^2(K_i)}$. 
Since $0 \leq \chi_j^{ms} \leq 1$, we have 
$0 \leq B \leq 1$ and $\vert \nabla B \vert^2 \leq C_\mathcal{T}^2 \Theta$. 
Using these facts together with $\nabla (B\nu)  = (\nabla B) \nu + B (\nabla \nu)$, we imply 
\begin{equation}
\begin{split}
\left\| \kappa^\frac{1}{2} \nabla (B\nu) \right\|_{L^2(K_i)}^2
& \leq C_\mathcal{T}^2 \|\nu\|_{s(K_i)}^2 + a_i(\nu,\nu) \\
& \leq \left( C_\mathcal{T}^2 + \lambda_{max} \right) \| \nu \|_{s(K_i)}^2 \\
\end{split}
\label{eq:inf-sup3}
\end{equation}
By taking the inf-sup constant
\begin{equation}
\beta = \dfrac{1 - C_\pi C_I \left( C_\mathcal{T}^2 + \lambda_{max} \right) \| \Theta\|_{L^\infty(K_i)}^\frac{1}{2} h}
{2C_\pi(1+C_I^2) \left( C_\mathcal{T}^2 + \lambda_{max} \right)}, 
\end{equation}
we prove the inf-sup condition \eqref{eq:inf-sup} and therefore the existence 
of $(v,\mu) \in V_{h,0}(K_i) \times V_{aux}^{(i)}$ in \eqref{eq:lemma2.0var}. 
It is then direct to check that the solution $v \in V_{h,0}(K_i)$ satisfies the desired properties. 
\end{proof}
\end{lemma}

We remark that, without loss of generality, we can assume $D \geq C_0^2(1+C_I^2)$.
We are now going to establish an estimate of the difference between 
the global multiscale basis functions and localized multiscale basis functions. 
We will see that the global multiscale basis functions have a decay property, 
and their values are small outside a suitably large oversampled domain. 
We will make use of a cutoff function in our proof. 
For each coarse block $K_i$ and $M > m$, 
the oversampling regions $K_{i,M}$ and $K_{i,m}$ 
define an outer neighborhood and an inner neighborhood respectively. 
We define $\chi_i^{M,m} \in \text{span}\{\chi_j^{ms}\}$ 
such that $0 \leq \chi_i^{M,m} \leq 1$ and 
\begin{equation}
\chi_i^{M,m} = 1 \text{ in } K_{i,m} \text { and }
\chi_i^{M,m} = 0 \text{ in } \Omega \setminus K_{i,M}.
\end{equation}
Moreover, we define the following DG norm for $w \in V$ 
on $K_{i,M} \setminus K_{i,m}$. 
\begin{equation}
\| w \|_{DG(K_{i,M} \setminus K_{i,m})}^2
= \sum_{K_k \subset K_{i,M} \setminus K_{i,m}} a_k(w,w) + 
\dfrac{\gamma}{h} \sum_{E \in \mathcal{E}^H(K_{i,M} \setminus K_{i,m})} 
\int_E \overline{\kappa} \llbracket w \rrbracket^2 d\sigma,
\end{equation}
where $\mathcal{E}^H(K_{i,M} \setminus K_{i,m})$ denotes the 
collection of all coarse grid edges in $\mathcal{E}^H$ which lie within 
in the interior of $K_{i,M} \setminus K_{i,m}$ and 
the boundary of $K_{i,M}$. 
We remark that the definition also applies to a region $\Omega \setminus K_{i,m}$ 
in the case when $M$ is sufficiently large. 

\begin{lemma}
\label{lemma3}
Assume the smallness criterion \eqref{eq:fine-small} on the fine mesh size $h$.
Suppose $m > 2$ is the number of coarse grid layers 
in the oversampled domain $K_{i,m}$ extended from 
the coarse grid block $K_i$. 
Let $\phi_j^{\left(i\right)} \in V_{aux}$ be a given auxiliary multiscale basis function. 
Let $\psi_j^{\left(i\right)} \in V_{glo}$ be the global multiscale basis function 
obtained from \eqref{eq:min1_glo}, and 
$\psi_{j,ms}^{\left(i\right)} \in V_h\left(K_{i,m}\right)$ be the localized multiscale basis function 
obtained from \eqref{eq:min1}.
Then we have
\begin{equation}
\| \psi_j^{\left(i\right)} - \psi_{j,ms}^{\left(i\right)} \|_a^2 \leq E \| \phi_j^{\left(i\right)} \|_{s(K_i)}^2,
\end{equation}
where $E = 40D^3 (1+\Lambda^{-1}) \left(1 + 6D^{-2 }\left(1 + \Lambda^{-\frac{1}{2}}\right)^{-1}\right)^{1-m}$.
\begin{proof}
By the variational formulations \eqref{eq:var1_glo} and \eqref{eq:var1}, we have 
\begin{equation}
a_{DG}\left(\psi_j^{(i)}-\psi_{j,ms}^{(i)},\psi\right) + s_i\left(\psi, \mu_j^{(i)} - \mu_{j,ms}^{(i)}\right) = 0 \text{ for all } \psi \in V_h(K_{i,m}). 
\label{eq:var1_diff}
\end{equation}
By Lemma~\ref{lemma2}, there exists $\widetilde{\phi}_j^{(i)} \in V_h$ such that 
\begin{equation}
\pi(\widetilde{\phi}_j^{(i)}) = \phi_j^{(i)}, \quad 
\| \widetilde{\phi}_j^{(i)} \|_a^2 \leq D \| \phi_j^{(i)} \|_{s(K_i)}^2, \quad 
\text{supp}\left(\widetilde{\phi}_j^{(i)}\right) \subseteq K_i.
\end{equation} 
We take $\eta = \psi_j^{(i)} - \widetilde{\phi}_j^{(i)} \in V_h$ and 
$\zeta =  \widetilde{\phi}_j^{(i)} - \psi_{j,ms}^{(i)} \in V_{h}(K_{i,m})$. 
By definition, we have $\pi(\eta) = \pi(\zeta) = 0$ and therefore 
$\eta, \zeta \in \widetilde{V}_h$. 
Again, by Lemma~\ref{lemma2}, there exists $\rho \in V_h$ such that 
\begin{equation}
\pi(\rho) = \pi(I_h(\chi_i^{m,m-1}\eta)), \quad 
\| \rho \|_a^2 \leq D \| \pi(I_h(\chi_i^{m,m-1}\eta)) \|_{s}^2, \quad 
\text{supp}\left(\rho\right) \subseteq K_{i,m} \setminus K_{i,m-1}.
\label{eq:def_rho}
\end{equation} 
Take $\tau = \rho - I_h(\chi_i^{m,m-1}\eta) \in V_h$. 
Again, $\pi(\tau) = 0$ and hence $\tau \in \widetilde{V}_h$. 
Taking $\psi = \tau - \zeta \in V_h(K_{i,m})$ in \eqref{eq:var1_diff} and 
making use of the fact $\tau - \zeta \in \widetilde{V}_h$, we have
\begin{equation}
a_{DG}\left(\psi_j^{(i)}-\psi_{j,ms}^{(i)},\tau-\zeta\right) = 0,
\end{equation} 
and therefore
\begin{equation}
\begin{split}
\left\| \psi_j^{\left(i\right)} - \psi_{j,ms}^{\left(i\right)} \right\|_a^2 
& = a_{DG}\left(\psi_j^{(i)}-\psi_{j,ms}^{(i)},\eta+\zeta\right) \\
& = a_{DG}\left(\psi_j^{(i)}-\psi_{j,ms}^{(i)},\eta+\tau\right) \\
& \leq \left\| \psi_j^{\left(i\right)} - \psi_{j,ms}^{\left(i\right)} \right\|_a 
\left\| \eta+\tau \right\|_a,
\end{split}
\end{equation}
which in turn implies 
\begin{equation}
\begin{split}
\left\| \psi_j^{\left(i\right)} - \psi_{j,ms}^{\left(i\right)} \right\|_a^2
& \leq \left\| \eta+\tau \right\|_a^2 \\
& = \left\| I_h((1-\chi_i^{m,m-1})\eta)+\rho \right\|_a^2 \\
& \leq 2\left( \left\| I_h((1-\chi_i^{m,m-1})\eta) \right\|_a^2 + \| \rho \|_a^2\right) \\
& \leq 2\left( C_0^2 \left\| I_h((1-\chi_i^{m,m-1})\eta) \right\|_{DG}^2 + \| \rho \|_a^2\right) \\
& \leq 2\left( 2C_0^2 (1+C_I^2) \left\| (1-\chi_i^{m,m-1})\eta \right\|_{DG}^2 + \| \rho \|_a^2\right).
\end{split}
\label{eq:lemma3.1}
\end{equation}
For the first term on the right hand side of \eqref{eq:lemma3.1}, 
by using $\nabla\left((1-\chi_i^{m,m-1})\eta\right) = -\nabla\chi_i^{m,m-1}\eta+(1-\chi_i^{m,m-1})\nabla\eta$ 
and $0 \leq 1-\chi_i^{m,m-1} \leq 1$, we have 
\begin{equation}
\begin{split}
\left\| (1-\chi_i^{m,m-1})\eta\right\|_a^2 \leq 2\left( 
\| \eta \|_{DG(\Omega \setminus K_{i,m-1})}^2 +
\left\|\eta\right\|_{s(\Omega \setminus K_{i,m-1})}^2\right).
\end{split}
\end{equation}
For the second term on the right hand side of \eqref{eq:lemma3.1}, using the definition of $\rho$ 
in \eqref{eq:def_rho} and $0 \leq 1-\chi_i^{m,m-1} \leq 1$, we obtain 
\begin{equation}
\| \rho \|_a^2 \leq D \| \pi (\chi_i^{m,m-1} \eta) \|_s^2 \leq D \| \chi_i^{m,m-1} \eta \|_s^2 \leq D \| \eta \|_{s(\Omega \setminus K_{i,m-1})}^2.
\end{equation}
Moreover, since $\eta \in \widetilde{V}_h$, by the spectral problem \eqref{eq:spectral_prob}, we have 
\begin{equation}
\| \eta \|_{s(\Omega \setminus K_{i,m-1})}^2 \leq \Lambda^{-1} \sum_{K_k \subset \Omega \setminus K_{i,m-1}} a_k(\eta,\eta).
\end{equation}
Combining all these estimates, we obtain 
\begin{equation}
\left\| \psi_j^{\left(i\right)} - \psi_{j,ms}^{\left(i\right)} \right\|_a^2
\leq 10D(1+\Lambda^{-1}) \| \eta \|_{DG(\Omega\setminus K_{i,m-1})}^2.
\label{eq:lemma3.2}
\end{equation}
Next, we will provide a recursive estimate for $\eta$ in the number of oversampling layers $m$. 
We take $\xi = 1- \chi_i^{m-1,m-2}$. Then 
$0 \leq \xi \leq 1$ and $\xi = 1$ in $\Omega \setminus K_{i,m-1}$. 
Using $\nabla(\xi^2 \eta) = \xi^2 \nabla \eta + 2\xi \eta \nabla \xi$, 
for every $K \in \mathcal{T}^H$, we have 
\begin{equation}
\int_K \kappa \nabla \eta \cdot \nabla(\xi^2 \eta) 
= \int_K \kappa \nabla \eta \cdot \left(\xi^2 \nabla \eta + 2\xi \eta \nabla \xi\right)
= \int_K \kappa \vert \nabla (\xi \eta) \vert^2 - \int_K \kappa \vert \nabla \xi \vert^2 \eta^2.
\end{equation}
In addition, using $\nabla(\xi\eta) = \xi \nabla \eta + \eta \nabla \xi$, 
for every $E \in \mathcal{E}^H$, we have 
\begin{equation}
\begin{split}
& -\int_E \{ \kappa \nabla \eta \cdot n_E \} \llbracket \xi^2 \eta \rrbracket 
-\int_E \{ \kappa \nabla (\xi^2 \eta) \cdot n_E \} \llbracket \eta \rrbracket + 
\dfrac{\gamma}{h} \int_E \overline{\kappa} \llbracket \eta \rrbracket \llbracket \xi^2 \eta \rrbracket  \\
& = -\int_E \{ \kappa \nabla \eta \cdot n_E \} \llbracket \xi^2 \eta \rrbracket 
-\int_E \left\{ \kappa (\xi^2 \nabla \eta + 2\xi \eta \nabla \xi) \cdot n_E \right\} \llbracket \eta \rrbracket + 
\dfrac{\gamma}{h} \int_E \overline{\kappa} \llbracket \eta \rrbracket \llbracket \xi^2 \eta \rrbracket \\
& = -2\left(\int_E \{ \kappa \xi \nabla \eta \cdot n_E \} \llbracket \xi\eta \rrbracket + 
\int_E \{ \kappa \eta \nabla \xi \cdot n_E \} \llbracket \xi \eta \rrbracket\right) + 
\dfrac{\gamma}{h} \int_E \overline{\kappa} \llbracket \xi \eta \rrbracket^2  \\
& = -2\int_E \{ \kappa \nabla (\xi \eta) \cdot n_E \} \llbracket \xi\eta \rrbracket +
\dfrac{\gamma}{h} \int_E \overline{\kappa} \llbracket \xi \eta \rrbracket^2.
\end{split}
\end{equation}
Summing over $K \in \mathcal{T}^H$ and $E \in \mathcal{E}^H$, we obtain 
\begin{equation}
\| \xi \eta \|_{a}^2 
\leq a_{DG}(\eta,\xi^2\eta) + \| \eta \|_{s(K_{i,m-1} \setminus K_{i,m-2})}^2,
\label{eq:lemma3.3}
\end{equation}
where we make use of the fact that $\nabla \xi = 0$ outside $K_{i,m-1} \setminus K_{i,m-2}$. 
We start with estimating the first term on the right hand side of \eqref{eq:lemma3.3}. 
For any coarse element $K_k \in \Omega \setminus K_{i,m-1}$, 
since $\xi = 1$ in $K_k$ and $\eta \in \widetilde{V}_h$, we have 
\begin{equation}
s\left(\xi^2 \eta, \phi_j^{(k)}\right) = s\left(\eta, \phi_j^{(k)}\right) = 0 \text{ for all } j = 1,2,\ldots,L_k.
\end{equation}
On the other hand, for any coarse element $K_k \in K_{i,m-2}$, 
since $\xi = 0$ in $K_k$, we have 
\begin{equation}
s\left(\xi^2 \eta, \phi_j^{(k)}\right) = 0 \text{ for all } j = 1,2,\ldots,L_k.
\end{equation}
Therefore, $\text{supp}(\pi(I_h(\xi^2 \eta))) \subset K_{i,m-1} \setminus K_{i,m-2}$. 
By Lemma~\ref{lemma2}, there exists $\sigma \in V_h$ such that 
\begin{equation}
\pi(\sigma) = \pi(I_h(\xi^2 \eta)), \quad 
\| \gamma \|_a^2 \leq D \| \pi (I_h(\xi^2 \eta)) \|_s^2, \quad 
\text{supp}(\sigma) \subset K_{i,m-1} \setminus K_{i,m-2}. 
\end{equation}
For any coarse element $K_k \subset K_{i,m-1} \setminus K_{i,m-2}$, 
since $0 \leq \xi \leq 1$ and $\pi(\eta) = 0$, we have 
\begin{equation}
\| \pi(I_h(\xi^2 \eta)) \|_{s(K_k)}^2 
\leq \| I_h(\xi^2 \eta) \|_{s(K_k)}^2 \leq \| I_h(\eta) \|_{s(K_k)}^2 
= \| \eta \|_{s(K_k)}^2 \leq \Lambda^{-1} a_k(\eta,\eta).
\end{equation}
Summing over $K_k \subset K_{i,m-1} \setminus K_{i,m-2}$, we obtain 
\begin{equation}
\| \pi(I_h(\xi^2 \eta)) \|_s^2 \leq \Lambda^{-1} \sum_{K_k \subset K_{i,m-1} \setminus K_{i,m-2}} a_k(\eta,\eta).
\end{equation}
We take $\theta = I_h(\xi^2 \eta) - \sigma$. 
Again, $\pi(\theta) = 0$ and $\theta \in \widetilde{V}_h$, which yields
\begin{equation}
a_{DG}\left(\psi_j^{(i)}, \theta\right) = 0.
\end{equation}
On the other hand, $\text{supp}(\theta) \subset \Omega \setminus K_{i,m-2}$ and $\text{supp}(\widetilde{\phi}_j^{(i)}) \subset K_i$. 
Since $\theta$ and $\widetilde{\phi}_j^{(i)}$ has disjoint supports, we have 
\begin{equation}
a_{DG}\left(\widetilde{\phi}_j^{(i)}, \theta\right) = 0.
\end{equation}
Therefore, we obtain 
\begin{equation}
a_{DG}(\eta,\theta) = a_{DG}\left(\psi_j^{(i)}-\widetilde{\phi}_j^{(i)}, \theta\right) = 0.
\end{equation}
Recall from the definition that $I_h(\xi^2 \eta) = \theta + \sigma$ and 
$\text{supp}(\sigma) \subset K_{i,m-1} \setminus K_{i,m-2}$.
Hence we have
\begin{equation}
\begin{split}
a_{DG}(\eta,I_h(\xi^2\eta)) 
& = a_{DG}(\eta,\sigma) \\
& \leq C_0\| \eta \|_{DG(K_{i,m-1}\setminus K_{i,m-2})} \| \sigma \|_a \\
& \leq C_0 D^\frac{1}{2} \| \eta \|_{DG(K_{i,m-1}\setminus K_{i,m-2})} \| \pi(I_h(\xi^2 \eta)) \|_s \\
& \leq D \Lambda^{-\frac{1}{2}} \| \eta \|_{DG(K_{i,m-1}\setminus K_{i,m-2})}^2.
\end{split}
\label{eq:lemma3.4a}
\end{equation}
On the other hand, making use of the fact that 
$\xi^2 = 0$ in $K_{i,m-2}$ and $\xi^2 = 1$ in $\Omega \setminus K_{i,m-1}$, 
we observe that $\xi^2 \eta = I_h(\xi^2 \eta)$ outside $K_{i,m-1} \setminus K_{i,m-2}$. 
Moreover, $\xi^2 \eta - I_h(\xi^2 \eta)$ is globally continuous. 
Thus, we obtain 
\begin{equation}
\begin{split}
a_{DG}(\eta,\xi^2\eta - I_h(\xi^2\eta)) 
& \leq C_0^2 \| \eta \|_{DG(K_{i,m-1}\setminus K_{i,m-2})} 
\left\|\xi^2\eta - I_h(\xi^2\eta) \right\|_{DG(K_{i,m-1}\setminus K_{i,m-2})}  \\
& \leq C_0^2 C_I \| \eta \|_{DG(K_{i,m-1}\setminus K_{i,m-2})} 
\left\|\xi^2\eta \right\|_{DG(K_{i,m-1}\setminus K_{i,m-2})} \\ 
& \leq \dfrac{D}{2} \left(\| \eta \|_{DG(K_{i,m-1}\setminus K_{i,m-2})}^2 + 
\left\|\xi^2\eta \right\|_{DG(K_{i,m-1}\setminus K_{i,m-2})}^2\right).
\end{split}
\label{eq:lemma3.4b}
\end{equation}
Again, using $\nabla(\xi^2 \eta) = \xi^2 \nabla \eta + 2\xi \eta \nabla \xi$, we have 
\begin{equation}
\left\|\xi^2\eta \right\|_{DG(K_{i,m-1}\setminus K_{i,m-2})}^2 \leq 
2 \| \eta \|_{DG(K_{i,m-1}\setminus K_{i,m-2})}^2 + 
8 \| \eta \|_{s(K_{i,m-1}\setminus K_{i,m-2})}^2.
\label{eq:lemma3.4c}
\end{equation}
Combining \eqref{eq:lemma3.3}, \eqref{eq:lemma3.4a}, \eqref{eq:lemma3.4b} and \eqref{eq:lemma3.4c}, we arrive at 
\begin{equation}
\| \xi \eta \|_{a}^2  \leq D \left(\left(\dfrac{3}{2} + \Lambda^{-\frac{1}{2}}\right) 
\| \eta \|_{DG(K_{i,m-1}\setminus K_{i,m-2})}^2 
+ 5\|  \eta \|_{s(K_{i,m-1}\setminus K_{i,m-2})}^2\right).
\end{equation}
Moreover, since 
$\pi(\eta) = 0$, we have 
\begin{equation}
\| \eta \|_{s(K_{i,m-1} \setminus K_{i,m-2})} \leq \Lambda^{-\frac{1}{2}} \| \eta \|_{DG(K_{i,m-1}\setminus K_{i,m-2})}, 
\end{equation}
which implies 
\begin{equation}
\| \xi \eta \|_{a}^2  \leq 6D \left(1 + \Lambda^{-\frac{1}{2}}\right) 
\|  \eta \|_{DG(K_{i,m-1}\setminus K_{i,m-2})}^2.
\end{equation}
By the equivalence of norms, we have 
\begin{equation}
\| \eta \|_{DG(\Omega\setminus K_{i,m-1})}^2 \leq C_0^2 \| \xi \eta \|_{a}^2 
\leq 6D^2 \left(1 + \Lambda^{-\frac{1}{2}}\right) \| \eta \|_{DG(K_{i,m-1}\setminus K_{i,m-2})}^2.
\end{equation}
We obtain the recurrence estimate 
\begin{equation}
\begin{split}
\| \eta \|_{DG(\Omega\setminus K_{i,m-2})}^2 
& = \| \eta \|_{DG(\Omega\setminus K_{i,m-1})}^2 + \| \eta \|_{DG(K_{i,m-1} \setminus K_{i,m-2})}^2 \\
& \geq \left(1 + 6D^{-2 }\left(1 + \Lambda^{-\frac{1}{2}}\right)^{-1}\right) \| \eta \|_{DG(\Omega\setminus K_{i,m-1})}^2.
\end{split}
\end{equation}
Inductively, we have 
\begin{equation}
\begin{split}
\| \eta \|_{DG(\Omega\setminus K_{i,m-1})}^2 
& \leq \left(1 + 6D^{-2 }\left(1 + \Lambda^{-\frac{1}{2}}\right)^{-1}\right)^{1-m} \| \eta \|_{DG(\Omega\setminus K_{i,1})}^2  \\
& \leq D \left(1 + 6D^{-2 }\left(1 + \Lambda^{-\frac{1}{2}}\right)^{-1}\right)^{1-m} \| \eta \|_{a}^2. 
\end{split}
\label{eq:lemma3.5}
\end{equation}
Combining \eqref{eq:lemma3.2} and \eqref{eq:lemma3.5}, we see that 
\begin{equation}
\left\| \psi_j^{(i)} - \psi_{j,ms}^{(i)} \right\|_a^2 \leq 10D^2 \left(1+\Lambda^{-1}\right)
\left(1 + 6D^{-2 }\left(1 + \Lambda^{-\frac{1}{2}}\right)^{-1}\right)^{1-m} \| \eta \|_{a}^2
\end{equation}
By the energy minimizing property of $\psi_j^{(i)}$, we have 
\begin{equation}
\| \eta \|_a \leq \| \psi_j^{(i)}\|_a + \| \widetilde{\phi}_j^{(i)} \|_a \leq 2 \| \widetilde{\phi}_j^{(i)} \|_a \leq 2D^\frac{1}{2} \| \phi_j^{(i)} \|_{s(K_i)}.
\end{equation}
We obtain the desired result. 
\end{proof}
\end{lemma}

Now, we are ready to establish our main theorem, 
which estimates the error between the solution $u_h$ 
and the multiscale solution $u_{ms}$.

\begin{theorem}
\label{theorem4}
Let $u_h \in V_h$ be the solution of \eqref{eq:sol_dg},
$u_{glo} \in V_{glo}$ be the solution of \eqref{eq:sol_glo}
with the global multiscale basis functions defined by \eqref{eq:min1_glo}, 
and $u_{ms} \in V_{ms}$ be the multiscale solution of \eqref{eq:sol_ms}
with the localized multiscale basis functions defined on an oversampled domain 
with $m > 2$ coarse grid layers by \eqref{eq:min1}.
Then we have 
\begin{equation}
\| u_h - u_{ms} \|_a \leq C \Lambda^{-\frac{1}{2}} \| \widetilde{\kappa}^{-\frac{1}{2}} f \|_{L^2(\Omega)} + 
Cm^d E^\frac{1}{2} \|u_{glo}\|_s,
\end{equation}
Moreover, if we let $k = O\left(\log\left(\dfrac{\kappa_1}{H}\right)\right)$ and 
replace the multiscale partition of unity $\{ \chi_j^{ms} \}$ by 
the bilinear partition of unity, we have
\begin{equation}
\| u_h - u_{ms} \|_a \leq C H \Lambda^{-\frac{1}{2}} \| \kappa^{-\frac{1}{2}} f \|_{L^2(\Omega)}.
\end{equation}
\begin{proof}
First, we write $u_{glo}$ in the linear combination of the basis $\{ \psi_k^{(j)} \}$ 
\begin{equation}
u_{glo} =  \sum_{i=1}^N \sum_{j=1}^{L_i} \alpha_j^{(i)} \psi_{j}^{(i)}.
\end{equation}
and define $\widehat{u}_{ms} \in V_{ms}$ by 
\begin{equation}
\widehat{u}_{ms} = \sum_{i=1}^N \sum_{j=1}^{L_i} \alpha_j^{(i)} \psi_{j,ms}^{(i)}.
\end{equation} 
From \eqref{eq:sol_dg} and \eqref{eq:sol_ms}, we obtain the Galerkin orthogonality 
\begin{equation}
a_{DG}(u_h - u_{ms}, w) = 0 \text{ for all } w \in V_{ms}, 
\end{equation}
which gives 
\begin{equation}
\| u_h - u_{ms} \|_{a} \leq \| u_h - \widehat{u}_{ms} \|_{a} \leq \| u_h - u_{glo} \|_{a} + \| u_{glo} - \widehat{u}_{ms} \|_{a}.
\label{eq:theorem4.1}
\end{equation}
Using Lemma~\ref{lemma3}, we see that
\begin{equation}
\begin{split}
\| u_{glo} - \widehat{u}_{ms} \|_{a}^2
& = \left\| \sum_{i=1}^N \sum_{j=1}^{L_i} \alpha_j^{(i)} (\psi_{j}^{(i)} - \psi_{j,ms}^{(i)}) \right\|_{a}^2 \\
& \leq Cm^d \sum_{i=1}^N \left\| \sum_{j=1}^{L_i} \alpha_j^{(i)} (\psi_{j}^{(i)} - \psi_{j,ms}^{(i)}) \right\|_{a}^2 \\
& \leq C m^d E \sum_{i=1}^N \left\| \sum_{j=1}^{L_i} \alpha_j^{(i)} \phi_j^{(i)} \right\|_s^2 \\
& = C m^d E \| u_{glo} \|_s^2,
\end{split}
\label{eq:theorem4.2}
\end{equation}
where the last equality follows from the orthogonality of the eigenfunctions in \eqref{eq:spectral_prob}. 
Using the estimates \eqref{eq:lemma1.0} and \eqref{eq:theorem4.2} in \eqref{eq:theorem4.1}, we have
\begin{equation}
\| u_h - u_{ms} \|_{a} 
\leq \Lambda^{-\frac{1}{2}} \| \widetilde{\kappa}^{-\frac{1}{2}} f \|_{L^2(\Omega)} + 
C m^{\frac{d}{2}} E^\frac{1}{2} \| u_{glo} \|_s.
\end{equation}
This completes the first part of the theorem. 
Next, we assume the partition of unity functions are bilinear, 
and we are going to estimate $\| u_{glo} \|_s$. 
Using the fact that $\vert \nabla \chi_k \vert = O(H^{-1})$, we have
\begin{equation}
\| u_{glo} \|_s^2 \leq CH^{-2} \kappa_1 \| u_{glo} \|_{L^2(\Omega)}^2.
\end{equation}
Then, by Poincar\'{e} inequality, we have
\begin{equation}
\| u_{glo} \|^2_{L^2(\Omega)} \leq C \kappa_0^{-1} \| u_{glo} \|_{a}^2.
\end{equation}
By taking $w = u_{glo} \in V_{glo}$ in \eqref{eq:sol_glo}, we obtain
\begin{equation}
\| u_{glo} \|_{a}^2 = (f, u_{glo})_{0,\Omega} \leq \| \widetilde{\kappa}^{-\frac{1}{2}} f \|_{L^2(\Omega)} \| u_{glo} \|_s.
\end{equation}
Combining these estimates, we have
\begin{equation}
\| u_{glo} \|_s \leq CH^{-2} \kappa_0^{-1} \kappa_1 \| \widetilde{\kappa}^{-\frac{1}{2}} f \|_{L^2(\Omega)}.
\end{equation}
To obtain our desired result, we need
\begin{equation}
H^{-2} \kappa_1 m^{\frac{d}{2}}E^\frac{1}{2} = O(1).
\end{equation}
Taking logarithm, we have
\begin{equation}
\log(H^{-2}) + \log(\kappa_1) + \dfrac{d}{2} \log(m) + \dfrac{1-m}{2} \log\left(1+\Lambda^{-\frac{1}{2}}\right) = O(1).
\end{equation}
Thus, taking $m = O\left(\log\left(\dfrac{\kappa_1}{H}\right)\right)$ 
completes the proof of the second result. 
\end{proof}
\end{theorem}

\section{Numerical results}
\label{sec:numerical}

In this section, we will present numerical examples with high contrast media 
to demonstrate the convergence of our proposed method 
with respect to the coarse mesh size $H$ and 
the number of oversampling layers $m$, 
and illustrate possible improvements in error robustness with respect to contrast by
employing the idea of constructing multiscale basis function by relaxation method introduced in \cite{chung2018constraint}.
Lastly, we examine the performance of applying the method to the wave equation. 
In all the experiments, the IPDG penalty parameter in \eqref{eq:dg_bilinear} 
is set to be $\gamma = 4$, so as to ensure the coercivity of the bilinear form $a_{DG}$.

\subsection{Experiment 1: flow problem}

\begin{figure}[ht!]
\centering
\includegraphics[width=0.6\linewidth]{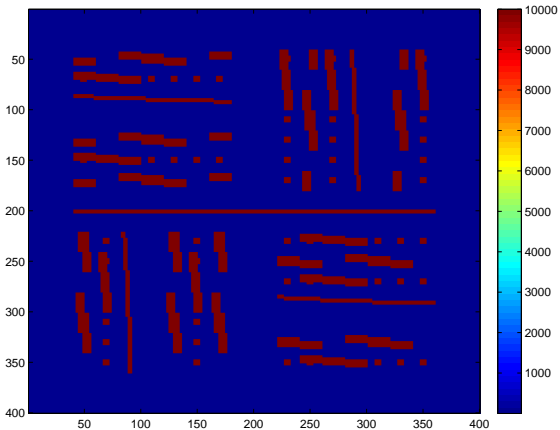}
\caption{The permeability field $\kappa$ for Experiment 1.}
\label{fig:medium1}
\end{figure}

In the first experiment, 
we consider a highly heterogeneous permeability field $\kappa$ 
in $\Omega = [0.1]^2$ as shown in Figure~\ref{fig:medium1}, 
with the background value is $\kappa = 1$ and 
the value in the channels and inclusions is $10^4$. 
and the resolution is $400\times400$, i.e.
$\kappa$ is piecewise constant on a fine grid with mesh size $h = 1/400$. 
The coarse mesh size varies from $H = 1/80$ to $H = 1/10$, 
and the number of oversampling layers varies from $m = 3$ to $m = 6$. 
In all these combinations, there are no more than $3$ high conductivity channels in a 
coarse block $K \in \mathcal{T}^H$. 
As a result, we have $3$ small eigenvalues in a local spectral problem \eqref{eq:spectral_prob}, 
and it suffices to use $3$ auxiliary basis functions per coarse block to construct the 
correspoding localized multiscale basis functions. 
The source function is taken as 
\begin{equation}
f(x,y) = 2\pi^2 \sin(\pi x) \sin(\pi y) \text{ for all } (x,y) \in \Omega.
\end{equation}
Table~\ref{tab:exp1.1} records the error when 
we take the number of oversampling layer to be approximately
$m \approx 4 \log(1/H)/ \log(1/10)$. 
The results show that the method provides optimal convergence in energy norm, 
which agrees with our theoretical finding in Section~\ref{sec:analysis},  
and the $L^2$ error converges with second order. 
Table~\ref{tab:exp1.2} records the error with various number of oversampling layers
and a fixed coarse mesh sizes $H = 1/40$.
It can be observed that increasing the number of oversampling layers 
improves the quality of approximations, but the decay in error is limited 
when the oversampling region is sufficiently large. 
This numerically verifies that the multiscale basis functions can indeed be localized. 

%


\begin{table}[ht!]
\centering
\begin{tabular}{cc|cc}
$m$ & $H$ & Energy error & $L^2$ error \\
\hline
4 & 1/10 & 7.4625\% & 0.7653\% \\
6 & 1/20 & 1.5392\% & 0.0625\% \\
7 & 1/40 & 0.7266\% & 0.0160\% \\
8 & 1/80 & 0.3433\% & 0.0035\% 
\end{tabular}
\caption{History of convergence with number of oversampling layers $m \approx 4 \log(1/H)/ \log(1/10)$ for Experiment 1.}
\label{tab:exp1.1}
\end{table}

\begin{table}[ht!]
\centering
\begin{tabular}{c|cc}
$m$ & Energy error & $L^2$ error \\
\hline
3 & 84.7517\% & 72.3079\% \\
4 & 19.0936\% & 3.6716\% \\
5 & 2.6687\% & 0.0720\% \\
6 & 0.7836\% & 0.0161\% \\
7 & 0.7266\% & 0.0160\% \\
8 & 0.7259\% & 0.0160\% 
\end{tabular}
\caption{Error table with different number of oversampling layers $m$ and a fixed coarse mesh size $H = 1/40$ for Experiment 1.}
\label{tab:exp1.2}
\end{table}

Next, we present the idea of the relaxed formulation of \eqref{eq:min1}.
Instead of using the method of Lagrange multiplier as in \eqref{eq:var1}, 
the $\phi$-orthogonality is imposed weakly by a penalty formulation.
The localized multiscale basis function $\psi_{j,{ms}}^{\left(i\right)} \in V_h\left(K_{i,m}\right)$ is defined as the solution of 
the following relaxed constrained energy minimization problem
\begin{equation}
\psi_{j,{ms}}^{\left(i\right)} = \argmin \left\{ a_{{DG}}\left(\psi, \psi\right) + s\left(\pi\left(\psi\right)-\phi_j^{\left(i\right)},\pi\left(\psi\right)-\phi_j^{\left(i\right)}\right) : 
\psi \in V_h\left(K_{i,m}\right)\right\}.
\label{eq:min2}
\end{equation}
The minimization problem \eqref{eq:min2} is equivalent to the following variational problem: 
find $\psi_{j,{ms}}^{\left(i\right)} \in V_h\left(K_{i,m}\right)$ such that
\begin{equation}
a_{{DG}}\left(\psi_{j,{ms}}^{\left(i\right)}, \psi\right) + s\left(\pi\left(\psi_{j,{ms}}^{\left(i\right)}\right), \pi\left(\psi\right)\right) = s\left(\phi_j^{\left(i\right)},\pi\left(\psi\right)\right) \text{ for all } \psi \in V_h\left(K_{i,m}\right).
\label{eq:var2}
\end{equation}
The construction of multiscale finite element space and coarse-scale model 
then follow \eqref{eq:msdg_space} and \eqref{eq:sol_ms} respectively. 
We compare the performance of the multiscale method with 
multiscale basis functions constructed by method of Lagrange multiplier \eqref{eq:var1} 
and the relaxation method \eqref{eq:var2} at different contrast values,  
where the coarse mesh size is taken as $H = 1/10$ 
and the number of oversampling layers as $m = 4$. 
In Table~\ref{tab:contrast}, we record the energy error and $L^2$ error 
with different contrast $\kappa_1$, 
where $\kappa_1 \gg 1$ is the value of $\kappa$ in the high conductivity channels.
It can be seen that the relaxation method is more robust with respect to contrast.

%

\begin{table}[ht!]
\centering
\begin{tabular}{c|cc|cc}
 & \multicolumn{2}{c|}{Lagrange multiplier} & \multicolumn{2}{c}{Relaxation} \\
\cline{2-5}
$\kappa_1$ & Energy error & $L^2$ error & Energy error & $L^2$ error \\
\hline
$10^4$ & 7.4625\% & 0.7653\% & 6.3757\% & 0.6395\% \\ 
$10^5$ & 12.6299\% & 1.6977\% & 6.3986\% & 0.6467\% \\ 
$10^6$ & 32.1465\% & 10.5146\% & 6.4020\% & 0.6478\% \\ 
$10^7$ & 64.1190\% & 41.8127\% & 6.4049\% & 0.6481\% \\ 
$10^8$ & 77.1229\% & 60.4947\% & 6.4301\% & 0.6503\% \\ 
\end{tabular}
\caption{Comparison of the method of Lagrange multiplier and the relaxation method with different contrast values for Experiment 1.}
\label{tab:contrast}
\end{table}

\subsection{Experiment 2: wave propagation}

\begin{figure}[ht!]
\centering
\includegraphics[width=0.6\linewidth]{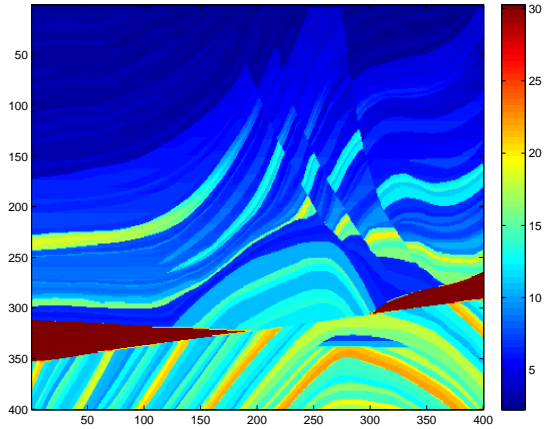}
\caption{Marmousi model for Experiment 2.}
\label{fig:marmousi}
\end{figure}

\begin{figure}[ht!]
\centering
\includegraphics[width=0.4\linewidth]{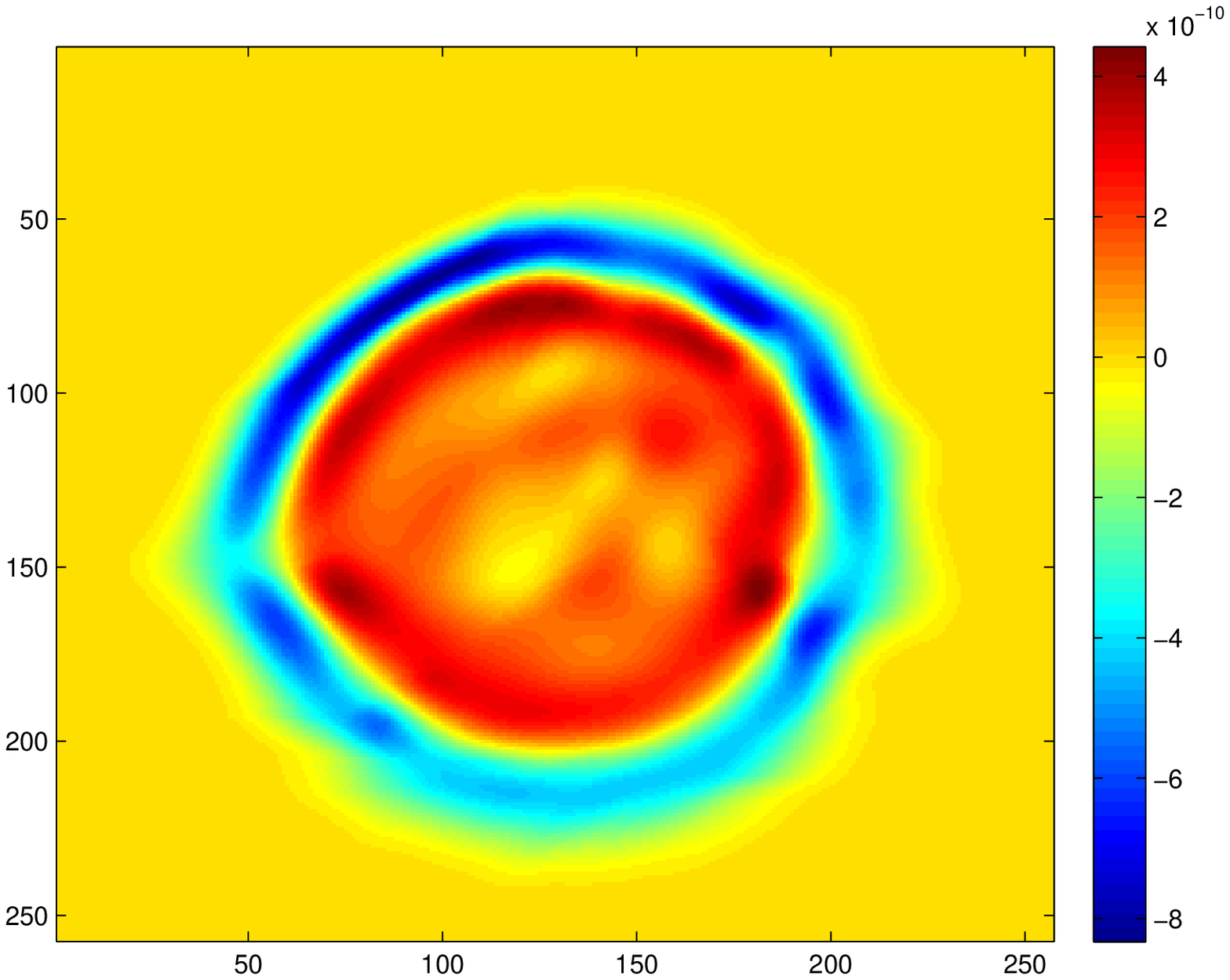}
\includegraphics[width=0.4\linewidth]{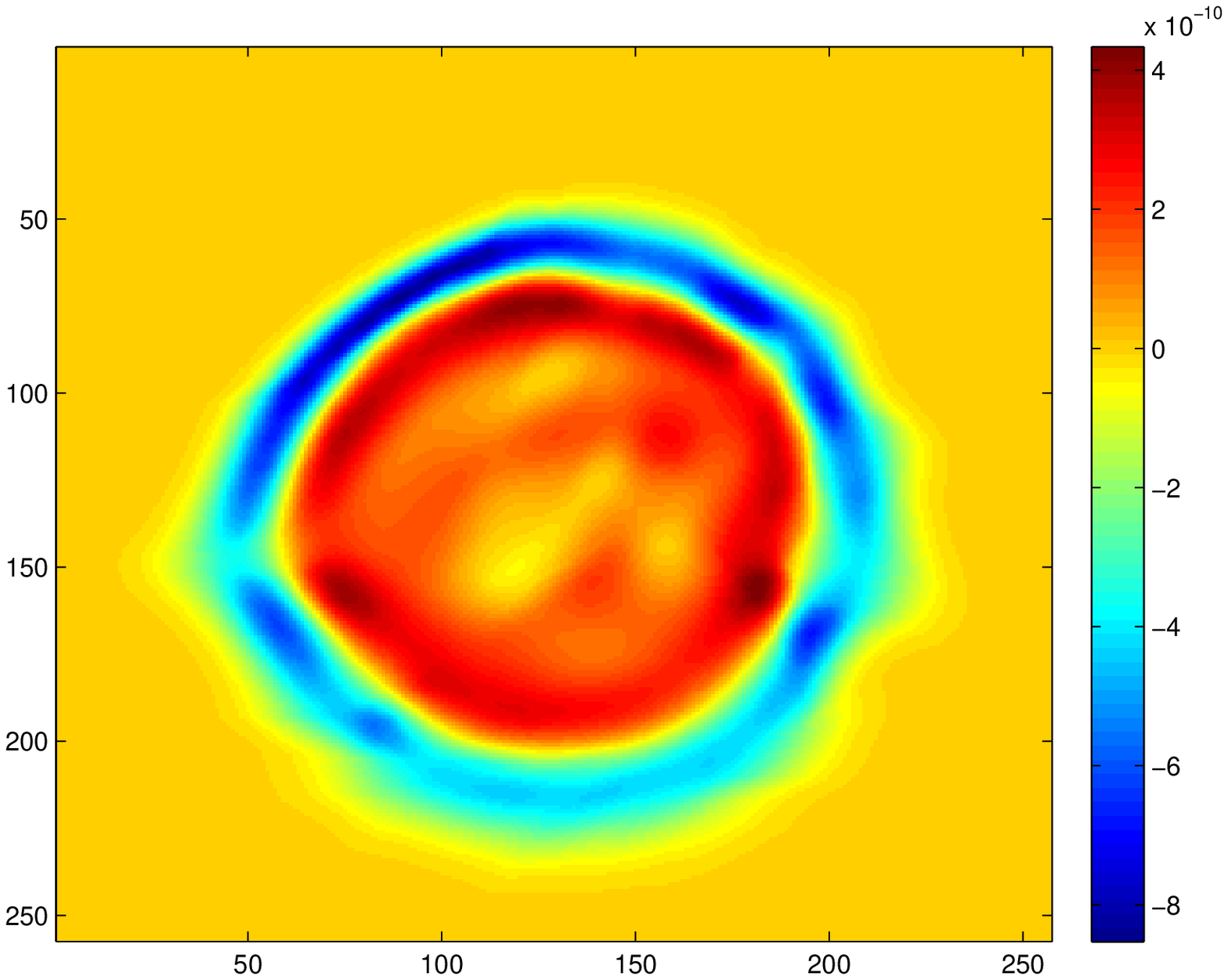}
\caption{Plots of numerical solution for Experiment 2. Fine solution (left) and multiscale solution (right).}
\label{fig:sol_wave}
\end{figure}

In the second experiment, we consider a simple wave equation in the space-time domain $\Omega_T = [0,T] \times \Omega$, 
where $[0,T]$ is the temporal domain and $\Omega$ is the spatial domain: 
\begin{equation}
\dfrac{\partial^2 u}{\partial t^2}  - \text{div}(\kappa \nabla u) = f,
\end{equation}
with homogeneous Dirichlet boundary condition $u = 0$ on $[0,T] \times \partial \Omega$. 
Here $\text{div}$ and $\nabla$ correspond to the divergence and gradient with respect to the spatial variable $x$, 
and $\kappa$ is the bulk modulus which is assumed to be stationary and highly oscillatory. 
Using second-order central difference for temporal discretization 
on a uniform grid points $t_k = k\Delta t$ and 
IPDG formulation on the space $V_h$, 
the fully-discrete DG finite element scheme reads: find $u^{n+1}_h \in V_h$ such that 
\begin{equation}
\dfrac{u_h^{n+1} - 2u_h^n + u_h^{n-1}}{\Delta t^2} + a_{DG}(u_h^n, w) = \int_\Omega f^n w \text{ for all } w \in V_h,
\end{equation}
where the superscript $k$ stands for the evaluation of a function at the time instant $t_k$, 
and the initial conditions $u^0_h, u^1_h$ are given. 
Next, we use the CEM method in Section~\ref{sec:method} 
to construct a multiscale finite element space $V_{ms}$. 
The coarse-scale model then reads: find $u^{n+1}_{ms} \in V_{ms}$ such that 
\begin{equation}
\dfrac{u_{ms}^{n+1} - 2u_{ms}^n + u_{ms}^{n-1}}{\Delta t^2} + a_{DG}(u_{ms}^n, w) = \int_\Omega f^n w \text{ for all } w \in V_{ms}.
\end{equation}
In this experiment, we take the bulk modulus on the spatial domain $\Omega = [0,1]^2$ 
as part of the Marmousi model as shown in Figure~\ref{fig:marmousi}. 
The fine mesh size is taken as $h = 1/256$. 
The coarse mesh size varies from $H = 1/32$ to $H = 1/8$, 
and the number of oversampling layers varies from $m = 4$ to $m = 6$. 
In all these combinations, we use $4$ auxiliary basis functions per coarse block to construct the 
correspoding localized multiscale basis functions. 
The source function $f$ is taken as the Ricker wavelet
\begin{equation}
f(t,x,y) = \dfrac{t-2/f_0}{4h^2} \exp\left( -\pi^2 f_0^2 (t-2/f_0)^2 \right) 
\exp\left(\dfrac{(x-0.5)^2+(y-0.5)^2}{4h^2}\right) \text{ for all } (t,x,y) \in \Omega_T.
\end{equation}
where the central frequency is chosen as $f_0 = 20$.
We iteratively solve for the numerical solution at the final time $T = 0.2$ 
with time step size $\Delta t = 10^{-4}$. 
Table~\ref{tab:error_wave} records the error of the final solution 
when we take the number of oversampling layer to be approximately
$m \approx 4 \log(1/H)/ \log(1/8)$. 
It can been observed that the method results in good accuracy and desired convergence in error. 
Figure~\ref{fig:sol_wave} depicts the numerical solutions 
by the fine-scale formulation and the coarse-scale formulation at the final time $T = 2$. 
The comparison suggests that the CEM method provides very good accuracy 
at a reduced computational expense.

\begin{table}[ht!]
\centering
\begin{tabular}{cc|cc}
$m$ & $H$ & Energy error & $L^2$ error \\
\hline
4 & 1/8 & 67.1831\% & 40.4876\% \\
5 & 1/16 & 24.1360\% & 9.7541\% \\
6 & 1/32 & 4.8839\% & 1.1559\%
\end{tabular}
\caption{History of convergence for Experiment 2.}
\label{tab:error_wave}
\end{table}

\section{Conclusion}
\label{sec:conclusion}
In this paper, we present CEM-GMsDGM, 
a local multiscale model reduction approach in the discontinuous Galerkin framework. 
The multiscale basis functions are defined in coarse oversampled regions by a constraint energy minimization problem, 
which are in general discontinuous on the coarse grid, and coupled by the IPDG formulation. 
Thanks to the definition of local spectral problems, 
the dimension of auxiliary space is minimal for sufficiently representing the high conductivity regions, 
and provides the most locally compressed multiscale space. 
In our analysis for the Darcy flow problem, we show that the method provides optimal convergence in the coarse mesh size, 
which is independent of the contrast, provided that the oversampling size is appropriately chosen. 
The convergence of the method for solving Darcy flow is theoretically analyzed and numerically verified. 
Numerical results for applying the method on a simple wave equation are also presented.

\bibliographystyle{plain}
\bibliography{references}

\begin{thebibliography}{10}

\bibitem{abdulle05}
A.~Abdulle.
\newblock On a priori error analysis of fully discrete heterogeneous multiscale
  fem.
\newblock {\em SIAM J. Multiscale Modeling and Simulation}, 4(2):447--459,
  2005.

\bibitem{buffa2006analysis}
A~Buffa, TJR Hughes, and G~Sangalli.
\newblock Analysis of a multiscale discontinuous {Galerkin} method for
  convection-diffusion problems.
\newblock {\em SIAM Journal on Numerical Analysis}, 44(4):1420--1440, 2006.

\bibitem{calo2011note}
Victor Calo, Yalchin Efendiev, and Juan Galvis.
\newblock A note on variational multiscale methods for high-contrast
  heterogeneous porous media flows with rough source terms.
\newblock {\em Advances in water resources}, 34(9):1177--1185, 2011.

\bibitem{ch03}
Z.~Chen and T.~Y. Hou.
\newblock A mixed multiscale finite element method for elliptic problems with
  oscillating coefficients.
\newblock {\em Mathematics of Computation}, 72(242):541--576, 2003.

\bibitem{cheung2018mc}
Siu~Wun Cheung, Eric~T Chung, Yalchin Efendiev, Wing~Tat Leung, and Maria
  Vasilyeva.
\newblock Constraint energy minimizing generalized multiscale finite element
  method for dual continuum model.
\newblock {\em arXiv preprint arXiv:1807.10955}, 2018.

\bibitem{cheung2019bayes}
Siu~Wun Cheung and Nilabja Guha.
\newblock Dynamic data-driven bayesian gmsfem.
\newblock {\em Journal of Computational and Applied Mathematics}, 353:72 -- 85,
  2019.

\bibitem{cgh09}
C.C. Chu, I.G. Graham, and T.~Hou.
\newblock A new multiscale finite element methods for high-contrast elliptic
  interface problem.
\newblock {\em Mathematics of Computation}, 79:1915--1955, 2010.

\bibitem{chung2016adaptive}
Eric Chung, Yalchin Efendiev, and Thomas~Y Hou.
\newblock Adaptive multiscale model reduction with generalized multiscale
  finite element methods.
\newblock {\em Journal of Computational Physics}, 320:69--95, 2016.

\bibitem{chung2018mixed}
Eric Chung, Yalchin Efendiev, and Wing~Tat Leung.
\newblock Constraint energy minimizing generalized multiscale finite element
  method in the mixed formulation.
\newblock {\em Computational Geosciences}, 22(3):677--693, 2018.

\bibitem{chung2015generalizedwave}
Eric~T Chung, Yalchin Efendiev, Richard~L Gibson~Jr, and Maria Vasilyeva.
\newblock A generalized multiscale finite element method for elastic wave
  propagation in fractured media.
\newblock {\em GEM-International Journal on Geomathematics}, pages 1--20, 2015.

\bibitem{chung2015residual}
Eric~T Chung, Yalchin Efendiev, and Wing~Tat Leung.
\newblock Residual-driven online generalized multiscale finite element methods.
\newblock {\em Journal of Computational Physics}, 302:176--190, 2015.

\bibitem{chung2017dg}
Eric~T. Chung, Yalchin Efendiev, and Wing~Tat Leung.
\newblock An online generalized multiscale discontinuous galerkin method
  (gmsdgm) for flows in heterogeneous media.
\newblock {\em Communications in Computational Physics}, 21(2):401–422, 2017.

\bibitem{chung2018constraint}
Eric~T Chung, Yalchin Efendiev, and Wing~Tat Leung.
\newblock Constraint energy minimizing generalized multiscale finite element
  method.
\newblock {\em Computer Methods in Applied Mechanics and Engineering},
  339:298--319, 2018.

\bibitem{chung2018dg}
E.T. Chung, Y.~Efendiev, and W.T. Leung.
\newblock An adaptive generalized multiscale discontinuous galerkin method for
  high-contrast flow problems.
\newblock {\em SIAM Multiscale Modeling and Simulation}, 16(3):1227--1257,
  2018.

\bibitem{chung2014adaptive}
ET~Chung, Y~Efendiev, and G~Li.
\newblock An adaptive {GMsFEM} for high-contrast flow problems.
\newblock {\em Journal of Computational Physics}, 273:54--76, 2014.

\bibitem{ee03}
W.~E and B.~Engquist.
\newblock Heterogeneous multiscale methods.
\newblock {\em Comm. Math. Sci.}, 1(1):87--132, 2003.

\bibitem{emz05}
W.~E, P.~Ming, and P.~Zhang.
\newblock Analysis of the heterogeneous multiscale method for elliptic
  homogenization problems.
\newblock {\em J. Amer. Math. Soc.}, 18(1):121--156, 2005.

\bibitem{egh12}
Y.~Efendiev, J.~Galvis, and T.~Hou.
\newblock Generalized multiscale finite element methods.
\newblock {\em Journal of Computational Physics}, 251:116--135, 2013.

\bibitem{eglmsMSDG}
Y~Efendiev, J~Galvis, R~Lazarov, M~Moon, and M~Sarkis.
\newblock Generalized multiscale finite element method. {Symmetric} interior
  penalty coupling.
\newblock {\em Journal of Computational Physics}, 255:1--15, 2013.

\bibitem{eglp13oversampling}
Y~Efendiev, J~Galvis, G~Li, and M~Presho.
\newblock Generalized multiscale finite element methods. {Oversampling}
  strategies.
\newblock {\em International Journal for Multiscale Computational Engineering,
  accepted}, 2013.

\bibitem{egw10}
Y.~Efendiev, J.~Galvis, and X.H. Wu.
\newblock Multiscale finite element methods for high-contrast problems using
  local spectral basis functions.
\newblock {\em Journal of Computational Physics}, 230:937--955, 2011.

\bibitem{eh09}
Y.~Efendiev and T.~Hou.
\newblock {\em Multiscale Finite Element Methods: Theory and Applications}.
\newblock Springer, 2009.

\bibitem{ehw99}
Y.~Efendiev, T.~Hou, and X.H. Wu.
\newblock Convergence of a nonconforming multiscale finite element method.
\newblock {\em SIAM J. Numer. Anal.}, 37:888--910, 2000.

\bibitem{efendiev2015spectral}
Yalchin Efendiev, Raytcho Lazarov, Minam Moon, and Ke~Shi.
\newblock A spectral multiscale hybridizable discontinuous {G}alerkin method
  for second order elliptic problems.
\newblock {\em Computer Methods in Applied Mechanics and Engineering},
  292:243--256, 2015.

\bibitem{efendiev2017bayes}
Yalchin Efendiev, Wing~Tat Leung, S.~W. Cheung, N.~Guha, V.~H. Hoang, and
  B.~Mallick.
\newblock Bayesian multiscale finite element methods. modeling missing subgrid
  information probabilistically.
\newblock {\em International Journal for Multiscale Computational Engineering},
  15(2):175--197, 2017.

\bibitem{elfverson2013dg}
D.~Elfverson, E.~Georgoulis, A.~Målqvist, and D.~Peterseim.
\newblock Convergence of a discontinuous galerkin multiscale method.
\newblock {\em SIAM Journal on Numerical Analysis}, 51(6):3351--3372, 2013.

\bibitem{hw97}
T.~Hou and X.H. Wu.
\newblock A multiscale finite element method for elliptic problems in composite
  materials and porous media.
\newblock {\em J. Comput. Phys.}, 134:169--189, 1997.

\bibitem{hou2017sparse}
Thomas~Y Hou and Pengchuan Zhang.
\newblock Sparse operator compression of higher-order elliptic operators with
  rough coefficients.
\newblock {\em Research in the Mathematical Sciences}, 4(1):24, 2017.

\bibitem{hfmq98}
T.J.R. Hughes, G.R. Feij\'oo, L.~Mazzei, and J.-B. Quincy.
\newblock The variational multiscale method - a paradigm for computational
  mechanics.
\newblock {\em Comput. Methods Appl. Mech Engrg.}, 127:3--24, 1998.

\bibitem{hughes2007variational}
TJR Hughes and G~Sangalli.
\newblock Variational multiscale analysis: the fine-scale {Green's} function,
  projection, optimization, localization, and stabilized methods.
\newblock {\em SIAM Journal on Numerical Analysis}, 45(2):539--557, 2007.

\bibitem{Iliev_MMS_11}
O.~Iliev, R.~Lazarov, and J.~Willems.
\newblock Variational multiscale finite element method for flows in highly
  porous media.
\newblock {\em Multiscale Model. Simul.}, 9(4):1350--1372, 2011.

\bibitem{maalqvist2014localization}
Axel M{\aa}lqvist and Daniel Peterseim.
\newblock Localization of elliptic multiscale problems.
\newblock {\em Mathematics of Computation}, 83(290):2583--2603, 2014.

\bibitem{owhadi2017multigrid}
Houman Owhadi.
\newblock Multigrid with rough coefficients and multiresolution operator
  decomposition from hierarchical information games.
\newblock {\em SIAM Review}, 59(1):99--149, 2017.

\bibitem{owhadi2014polyharmonic}
Houman Owhadi, Lei Zhang, and Leonid Berlyand.
\newblock Polyharmonic homogenization, rough polyharmonic splines and sparse
  super-localization.
\newblock {\em ESAIM: Mathematical Modelling and Numerical Analysis},
  48(2):517--552, 2014.

\bibitem{papanicolau1978asymptotic}
G~Papanicolau, A~Bensoussan, and J-L Lions.
\newblock {\em Asymptotic analysis for periodic structures}.
\newblock Elsevier, 1978.

\bibitem{park2019mc}
Jun Sur~Richard Park, Siu~Wun Cheung, Tina Mai, and Viet~Ha Hoang.
\newblock Multiscale simulations for upscaled multi-continuum flows.
\newblock {\em arXiv preprint arXiv:1909.04722}, 2019.

\bibitem{riviere2008discontinuous}
B{\'e}atrice Rivi{\`e}re.
\newblock {\em Discontinuous Galerkin methods for solving elliptic and
  parabolic equations: theory and implementation}.
\newblock Society for Industrial and Applied Mathematics, 2008.

\bibitem{wang2020vug}
Min Wang, Siu~Wun Cheung, Eric~T. Chung, Maria Vasilyeva, and Yuhe Wang.
\newblock Generalized multiscale multicontinuum model for fractured vuggy
  carbonate reservoirs.
\newblock {\em Journal of Computational and Applied Mathematics}, 366:112370,
  2020.

\bibitem{weh02}
X.H. Wu, Y.~Efendiev, and T.Y. Hou.
\newblock Analysis of upscaling absolute permeability.
\newblock {\em Discrete and Continuous Dynamical Systems, Series B.},
  2:158--204, 2002.

\end{thebibliography}

\end{document}